\documentclass[reqno,11pt]{amsart}
\usepackage{amscd,amssymb,verbatim}

\numberwithin{equation}{section}
\theoremstyle{plain}

\newcommand\alp{\alpha}         
\newcommand\bet{\beta}
         \newcommand\Gam{\Gamma}
         
\newcommand\eps{\varepsilon}

\newcommand\iot{\iota}

                \newcommand\Lam{\Lambda}
         \newcommand\Sig{\Sigma}

\newcommand\ome{\omega}         

\newcommand\calE{{\mathcal{E}}}

\newcommand\calH{{\mathcal{H}}}

\newcommand\calL{{\mathcal{L}}}

            \newcommand\bfD{{\mathbf D}}

\newcommand\bfu{{\mathbf u}}

\newcommand\RR{\mathbb{R}}

\newcommand\ZZ{\mathbb{Z}}

\newcommand\CC{\mathbb{C}}

 \newcommand\grg{{\mathfrak{g}}}

\newcommand\sdp{\times \hskip -0.3em {\raise 0.3ex
\hbox{$\scriptscriptstyle |$}}} % semidirect product

\newcommand\Ad{\operatorname{Ad\, }}

\newcommand\End{\operatorname{End\,}}

\newcommand\ind{\operatorname{ind}}

\newcommand\Ker{\operatorname{Ker}}

\newcommand\Lie{\operatorname{Lie}}

\newcommand\spin{\operatorname{spin}}

\renewcommand\oe{{\bar{e}}}

\newcommand\tilc{{\tilde{c}}}

\newcommand\tilD{{\tilde{D}}}
\newcommand\tilE{{\tilde{E}}}

\newcommand\tilv{{\tilde{v}}}

\newcommand\tilW{{\tilde{W}}}

\renewcommand{\>}{\rangle}
\newcommand{\<}{\langle}

\theoremstyle{plain}
\newtheorem{Thm}[subsection]{Theorem}
\newtheorem{Cor}[subsection]{Corollary}
\newtheorem{Lem}[subsection]{Lemma}
\newtheorem{Prop}[subsection]{Proposition}
\newtheorem{Conjec}[subsection]{Conjecture}

\theoremstyle{definition}
\newtheorem{Def}[subsection]{Definition}%\renewcommand{\theDef}{\thesection.\arabic{Def}}

\theoremstyle{remark}

\newtheorem{Rem}[subsection]{Remark}%\renewcommand{\theRem}{\thesection.\arabic{Rem}}

\newif\ifShowLabels
\ShowLabelstrue
\newdimen\theight
\def\TeXref#1{%
        \leavevmode\vadjust{\setbox0=\hbox{{\tt
                \quad\quad  {\small \textrm #1}}}%
        \theight=\ht0
        \advance\theight by \lineskip
        \kern -\theight \vbox to
        \theight{\rightline{\rlap{\box0}}%
        \vss}%
        }}%

\newcommand{\refs}[1]{Section ~\ref{S:#1}}
\newcommand{\refss}[1]{Subsection ~\ref{SS:#1}}
\newcommand{\reft}[1]{Theorem ~\ref{T:#1}}
\newcommand{\refl}[1]{Lemma ~\ref{L:#1}}
\newcommand{\refp}[1]{Proposition ~\ref{P:#1}}

\newcommand{\refd}[1]{Definition ~\ref{D:#1}}
\newcommand{\refr}[1]{Remark ~\ref{R:#1}}
\newcommand{\refe}[1]{\eqref{E:#1}}

\newenvironment{thm}[1]%
        { \begin{Thm} \label{T:#1}  \ifShowLabels \TeXref{T:#1} \fi }%
        { \end{Thm} }

\renewcommand{\th}[1]{\begin{thm}{#1} \sl }
\renewcommand{\eth}{\end{thm} }

\newenvironment{lemma}[1]%
        { \begin{Lem} \label{L:#1}  \ifShowLabels \TeXref{L:#1} \fi }%
        { \end{Lem} }
\newcommand{\lem}[1]{\begin{lemma}{#1} \sl}
\newcommand{\elem}{\end{lemma}}

\newenvironment{propos}[1]%
        { \begin{Prop} \label{P:#1}  \ifShowLabels \TeXref{P:#1} \fi }%
        { \end{Prop} }
\newcommand{\prop}[1]{\begin{propos}{#1}\sl }
\newcommand{\eprop}{\end{propos}}

\newenvironment{corol}[1]%
        { \begin{Cor} \label{C:#1}  \ifShowLabels \TeXref{C:#1} \fi }%
        { \end{Cor} }
\newcommand{\cor}[1]{\begin{corol}{#1} \sl }
\newcommand{\ecor}{\end{corol}}

\newenvironment{conjec}[1]%
        { \begin{Conjec} \label{Conj:#1}  \ifShowLabels \TeXref{C:#1} \fi }%
        { \end{Conjec} }
\newcommand{\conj}[1]{\begin{conjec}{#1} \sl }
\newcommand{\econj}{\end{conjec}}

\newenvironment{remark}[1]%
        { \begin{Rem} \label{R:#1}  \ifShowLabels \TeXref{R:#1} \fi }%
        { \end{Rem} }
\newcommand{\rem}[1]{\begin{remark}{#1}}
\newcommand{\erem}{\end{remark}}

\newcommand{\eq}[1]%
        { \ifShowLabels\newline \TeXref{E:#1} \fi
           \begin{equation} \label{E:#1} }
\newcommand{\eeq}{\end{equation}}

\newcommand{\prf}{ \begin{proof} }
\newcommand{\eprf}{ \end{proof} }
\newcommand{\Label}[1]{\label{#1}  \ifShowLabels \TeXref{#1} \fi }

\ShowLabelsfalse% comment this out if labels should be printed

\renewcommand{\d}{\text{\( \partial\)}}

\renewcommand{\b}{\bullet}

\newcommand{\n}{\nabla}

\newcommand{\E}{\calE}

\renewcommand{\L}{\calL}

\newcommand{\g}{{\Gam}}
\newcommand{\gc}{{\Gam(M,C(M))}}
\newcommand{\gme}{{\Gam(M,\E)}}

\newcommand{\nLC}{\n^{\text{LC}}}
\newcommand{\wE}{\tilde\calE}

\renewcommand{\v}{\mathbf{v}}\renewcommand{\u}{\mathbf{u}}
\newcommand{\tv}{\mathbf{\tilv}}
\renewcommand{\i}{\sqrt{-1}\, }

\newcommand{\Gtc}{\Gam_{tc}^\infty(M,\E)}

\begin{document}

%\hfill{Preliminary version: \today}
\title{The index theory on non-compact manifolds with  proper group action}
\author{Maxim Braverman}
\thanks{Supported in part by the NSF grant DMS-1005888.}
\address{Department of Mathematics,
        Northeastern University,
        Boston, MA 02115,
        USA
         }

%\date{\today}
\begin{abstract}
We construct a regularized index of a generalized Dirac operator on a complete Riemannian manifold endowed with a proper action of a unimodular Lie group. We show that the index is preserved by a certain class of non-compact cobordisms and prove a gluing formula for the regularized index. The results of this paper generalize our previous construction of index for compact group action and the recent paper of Mathai and Hochs who studied  the case of a Hamiltonian action on a symplectic manifold. As an application of the cobordism invariance of the index we give an affirmative answer to a question of Mathai and Hochs about the independence of the Mathai-Hochs quantization of the metric, connection and other choices.  
\end{abstract}
\maketitle

%------------------------------------------------------
%------------------------------------------------------
\section{Introduction}\label{S:introduction}

Paradan \cite{Paradan03}  introduces a regularized topological index of a Dirac-type operator on a non-compact manifold $M$ endowed with an action of a compact group $G$. In  \cite{Br-index} (see also \cite{BrCano} for a review) we constructed an analytic counterpart of this index and proved that the two indexes coincide (see also \cite{MaZhang_TrIndex12}). The regularized index depends on an additional data, namely an equivariant map from $M$ to the Lie algebra of $G$, called the {\em taming map}.  The regularized index was used in  \cite{MaZhang-noncompact}  as a method to  prove a conjecture of Vergne \cite{Vergne07}. The method of \cite{Br-index} was also used in \cite{BrBackground} to construct  a regularized Dolbeault cohomology of a non-compact $G$-manifold.

Mathai and Zhang \cite{MathaiZhang10} defined an index of a Dirac operator on a manifold endowed with a proper cocompact action of a non-compact group. Mathai and Hochs  \cite{MathaiHochs} considered a Hamiltonian action of a non-compact group $G$ on a non-compact symplectic manifold $M$. They constructed a regularized  index of the $\spin^c$-Dirac operator without assuming that the action is cocompact and showed that this index has the Guillemin-Sternberg ``quantization commutes with reduction" property. The regularization of the index used in \cite{MathaiHochs} is very similar to the one introduced in \cite{Br-index}, with the moment map playing the role of the taming map. In this sense, the Mathai-Hochs construction can be viewed  as a combination of ideas from \cite{MathaiZhang10} and \cite{Br-index}. 

In this note  we combine the methods of \cite{Br-index} and \cite{MathaiHochs} to construct an analytical index of a generalized Dirac operator on a  non-compact manifold endowed with a proper action of a non-compact Lie group. In the case of a Hamiltonian group action on a symplectic manifold our index coincides with the construction of \cite{MathaiHochs}. We show that our index is invariant under a certain class of non-compact cobordisms. We also prove a gluing formula for this index. From the cobordism invariance of the index we immediately conclude that the index is independent of the metric, the connection and other data used in its definition. That gives an affirmative answer to a question posed by Mathai and Hochs, cf. Remarks~3.8 and  6.2 of \cite{MathaiHochs}.

%------------------
\subsection{The construction of the index}\label{SS:Iconstruction}
Suppose $M$ is a complete Riemannian manifold on which a non-compact unimodular Lie group $G$ acts by isometries. Let $\E=\E^+\oplus\E^-$ be a $G$-equivariant $\ZZ_2$-graded self-adjoint Clifford module over $M$. We refer to the pair $(\E,\v)$ as a {\em tamed Clifford module}.

Consider a Dirac operator $D^\pm:\Gam(M,\E^\pm)\to \Gam(M,\E^\mp)$ associated to a Clifford connection on $\E$. 

Following \cite{MathaiZhang10,MathaiHochs} we consider a smooth {\em cutoff function} $\chi:M\to [0,\infty)$, whose support intersects all $G$-orbits in compact sets and which satisfies $\int_G  \chi(g\cdot x)^2\, dg= 1$. Since the group $G$ is unimodular, such a function always exists by \cite[Ch. VII,\S2.4]{Bourbaki-integration7-8}.

A section  $s\in \Gam(M,\E)$ is called {\em transversally compactly supported}\/ if the support of $s$ is cocompact. We denote by $\Gtc$ the space of smooth transversally compactly supported sections of $\E$ and by $\Gtc^G$ the subspace of  $G$-invariant elements of $\Gtc$.
Consider the operator 
\[
	D_{\chi}:\, \chi\,\Gtc^G\ \to \ \chi\, \Gtc^G, \qquad D_{\chi}(\chi s)\ := \ \chi\,Ds, \qquad \text{for} \ \ s\in \Gtc.
\]
It is shown in \cite{MathaiZhang10} that if the quotient space $M/G$ is compact than the operator $D_\chi$ is Fredholm and its index is independent of the choice of the function $\chi$. Thus one can define the regularized index of $D$ by 
\(
	\ind_G D\ := \ \ind D_\chi
\)
for any cutoff function $\chi$. If $M/G$ is not compact then $D_\chi$ is not necessarily Fredholm and a regularization is needed to define its index. 

Let $\v:M\to\grg=\Lie G$ be a $G$-equivariant map, such that the induced vector field $v$ on $M$ does not vanish outside of a cocompact subset of $M$.  We call $\v$ a {\em taming} map, and we refer to the pair $(\E,\v)$ as a {\em tamed Clifford module}.

Let $f:M\to[0,\infty)$ be a $G$-invariant function which increases fast enough at infinity
(see \refss{rescaling} for the precise condition on $f$). We
consider the {\em deformed Dirac operator} 
\[
	D_{\chi,fv} \ := \ D_\chi+{\i}c(fv):\, \chi\,\Gtc^G\ \to \ \chi\, \Gtc^G,
\]	 
where $c:TM\simeq T^*M\to\End\E$ is the Clifford module structure
on $\E$.  Our principal result is  \reft{finite}, which states that the deformed Dirac operator is Fredholm and its index $\ind_GD_{\chi,fv}$ is independent  of
the choice of the functions $\chi, \ f$ and of the Clifford connection on $\E$, used in the definition of $D$.  We denote this index by $ \ind_G(\E,\v)$ and call it  the {\em (analytic) index of $(\E,\v)$}.

%--------------------
\subsection{The cobordism invariance}\label{SS:Icobordim}
In \refs{cobord}, we introduce the notion of a cobordism between tamed Clifford modules. Roughly speaking, this is a usual cobordism, which carries a taming map. Our notion of cobordism is very close to the notion of non-compact cobordism developed by V.~Ginzburg, V.~Guillemin and Y.~Karshon \cite{GGK96,Karshon98,GGK-book}. We prove, that {\em the index is preserved by a cobordism}. 

%--------------
\subsection{The gluing formula}\label{SS:Igluing}
Suppose $\Sig\subset M$ is a cocompact $G$-invariant hypersurface,
such that the vector field $v$ does not vanish anywhere on $\Sig$.
We endow the open manifold $M\backslash{\Sig}$ with a complete
Riemannian metric and we denote by $(\E_\Sig,\v_\Sig)$ the induced
tamed Clifford module on $M\backslash{\Sig}$. In \refs{gluing}, we
prove that {\em  the tamed Clifford modules $(\E_\Sig,\v_\Sig)$
and $(\E,\v)$ are cobordant. In particular, they have the same
index}. We refer to this result as the {\em gluing formula}. 

The gluing formula takes especially nice form if $\Sig$ divides
$M$ into 2 disjoint manifolds $M_1$ and $M_2$. Let $(\E_1,\v_1)$
and $(\E_2,\v_2)$ be the restrictions of $(\E_\Sig,\v_\Sig)$ to
$M_1$ and $M_2$, respectively. Then the gluing formula implies
\[
    \ind_G(\E,\v) \ = \ \ind_G(\E_1,\v_1) \ +  \ \ind_G(\E_2,\v_2).
\]
In other words, {\em the index is additive}. 
%-----------------------------------------------
%-----------------------------------------------
\section{Equivariant Index for a unimodular group action on non-compact manifolds}\label{S:defofindex}

In this section we introduce our main objects of study: tamed
non-compact manifolds, tamed Clifford modules, and the (analytic)
equivariant index of such modules.

Throughout the paper $(M,g^M)$ is a complete Riemannian manifold without boundary. 

\subsection{Clifford module and Dirac operator}\Label{SS:dirac}
Let $C(M)$ denote the Clifford bundle of $M$ (cf. \cite[\S3.3]{BeGeVe}),
i.e., a vector bundle, whose fiber at every point $x\in M$ is
isomorphic to the Clifford algebra $C(T^*_xM)$ of the cotangent
space.

A {\em ($\ZZ_2$-graded self-adjoint) Clifford module \/} on $M$ is a $\ZZ_2$-graded Hermitian  vector bundle $\E=\E^+\oplus\E^-$ over $M$ endowed with a graded action
\[
        (a,s) \ \mapsto \ c(a)s, \quad \mbox{where} \quad
                        a\in \gc, \ s\in \gme,
\]
of the bundle $C(M)$ such that the operator $c(v):\E_x\to\E_x$ is skew-adjoint, for all $x\in M$ and $v\in T_x^*M$.

A {\em Clifford connection}  on $\E$ is a Hermitian connection
$\n^\E$, which preserves the subbundles $\E^\pm$ and
\[
    [\n^\E_X,c(a)] \ = \ c(\nLC_X a), \quad
                \mbox{for any} \quad  a\in \gc, \ X\in\g(M,TM),
\]
where $\nLC_X$ is the Levi-Civita covariant derivative on $C(M)$
associated with the Riemannian metric on $M$.

The {\em Dirac operator \/} $D:\gme\to\gme$ associated to a
Clifford connection $\n^\E$ is defined by the following
composition
\[
  \begin{CD}
        \gme @>\n^\E>> \g(M,T^*M\otimes \E) @>c>> \gme.
  \end{CD}
\]
In local coordinates, this operator may be written as
$D=\sum\,c(dx^i)\,\n^\E_{\d_i}$. Note that $D$ sends even sections
to odd sections and vice versa: $D:\, \Gam(M,\E^\pm)\to
\Gam(M,\E^\mp)$.

Consider the $L^2$-scalar product on the space of sections $\gme$
defined by the Riemannian metric on $M$ and the Hermitian
structure on $\E$. By \cite[Proposition~3.44]{BeGeVe}, the Dirac
operator associated to a Clifford connection $\n^\E$ is formally
self-adjoint with respect to this scalar product. Moreover, it is
essentially self-adjoint with the initial domain smooth, compactly
supported sections, cf. \cite{Chernoff73}, \cite[Th.~1.17]{GromLaw83}.

If $W$ is a manifold with boundary, then by a Clifford module over it we will understand
a smooth vector bundle $\E$ over $W$, whose restriction to the set of the interior points
$W^{\text{int}}$ of $W$ has a structure of a Clifford module over $W^{\text{int}}$.
(Usually we require some additional structure of $\E$ near the boundary of $W$, but
we will formulate these requirements when we need them.)

%--------------------------------
\subsection{Group action. The index.}\Label{SS:G}
Let $G$ be a unimodular  Lie group and suppose that there is a proper action of $G$ on $M$ by isometries. Assume that there is given a lift of this action to $\E$, which
preserves the grading, the connection and the Hermitian metric on $\E$. Then the Dirac operator $D$ commutes with the action of $G$.
Hence, $\Ker D= \Ker D^+\oplus\Ker{}D^-$ is a $G$-invariant subspace of  $\Gam(M,\E)$ . Let 
\[
	(\Ker{}D^\pm)^G\ \subset\ \Ker D^\pm
\]
denote the space of $G$-invariant elements of $\Ker D^\pm$.

If $M$ is compact, then the spaces  $\Ker D^\pm$ and, hence, $(\Ker D^\pm)^G$  are finite dimensional.  In this situation we define an {\em equivariant index} of $D$ by 
\begin{equation}\label{E:compactindex}
        \ind_G(D) \ = \ \dim (\Ker{}D^+)^G \ - \ \dim (\Ker{}D^-)^G.
\end{equation}
The index  depends only  on $M$ and the equivariant Clifford module $\E=\E^+\oplus\E^-$ and does not
depend on the choice of the connection $\n^\E$ and the metric $h^\E$.  We set 
\[
	\ind_G(\E)\ := \ \ind_G(D),
\]
and refer to it as the {\em index} of $\E$.

Our goal is to define  an analogue of \refe{compactindex} for a $G$-equivariant Clifford module over a complete {\em non-compact \/} manifold. If the group $G$ is compact, then this was done on \cite{Br-index}.  It the group $G$ is not compact, but $M$ is a symplectic manifold and the action of $G$ on $M$ is Hamiltonian, the index was constructed in \cite{MathaiHochs}. We now combine the ideas of \cite{Br-index} and \cite{MathaiHochs} to construct the index for a general complete $G$-manifold. 

%---------------------------------------------
\subsection{Cutoff function along the orbits}\label{SS:cutoff} 
Let $dg$ denote the Haar measure on $G$. Following \cite{MathaiZhang10,MathaiHochs} we make the following definition 
\begin{Def}\label{D:cutoff}
 A smooth  function $\chi:M\to [0,\infty)$, whose support intersects all $G$-orbits in compact sets and which satisfies 
\begin{equation}\label{E:cutoff}
	\int_G  \chi(g\cdot x)^2\, dg\ = \ 1
\end{equation}
is called a {\em cutoff function} on $M$. 
\end{Def}
Note that a cutoff function always exists by \cite[Ch. VII,\S2.4]{Bourbaki-integration7-8}.

Let $M/G$ denote the space of $G$-orbits in $M$ and let $q:M\to M/G$ denote the quotient map.  

%%%%
\begin{Def}\label{D:trcompactset}
A subset $V\subset M$ is called {\em cocompact} if $q(V)\subset M/G$ is compact.
\end{Def}
\begin{Def}\label{D:trcompactsupport}
A section $s\in \Gam(M,\E)$ is called {\em transversally compactly supported}\/ if the support of $s$ is cocompact. We denote by $\Gtc$ the space of smooth transversally compactly supported sections of $\E$ and by $\Gtc^G$ the subspace  of $G$-invariant elements in $\Gtc$.
\end{Def}
Notice that if $s\in \Gtc$ then $\chi{}s$ is a smooth compactly supported section of $\E$. 

Define the operator 
\[
	D_{\chi}:\, \chi\,\Gtc^G\ \to \ \chi\, \Gtc^G
\]
by
\begin{equation}\label{E:Dphi}
	D_{\chi}(\chi s)\ := \ \chi\,Ds, \qquad s\in \Gtc.
\end{equation}

It is shown in \cite{MathaiZhang10} that if the quotient space $M/G$ is compact than the operator $D_\chi$ is Fredholm and its index is independent of the choice of the function $\chi$. Thus one can define the regularized index of $D$ by 
\[
	\ind_G D\ := \ \ind D_\chi
\]
for any cutoff function $\chi$. As before, $\ind_G(D)$ depends only on $M$ and the equivariant Clifford module $\E=\E^+\oplus\E^-$ and we set  $\ind_G(\E)\ := \ \ind_G(D)$.

%------------------------------------
\subsection{A tamed non-compact manifold} \label{SS:tamed}
When $M/G$ is not compact to define the index we need and additional structure on $M$. This structure is given by an equivariant map $\v:M\to\grg$, where $\grg$ denotes the Lie algebra of $G$ and $G$ acts on it by the adjoint representation. Such a map induces a vector field $v$ on $M$ defined by
\begin{equation}\label{E:v}
    v(x) \ := \ \frac{d}{dt}\Big|_{t=0}\, \exp{(t\v(x))}\cdot x.
\end{equation}

The following definition extends Defenition~2.4 of \cite{Br-index}:

\begin{Def}\label{D:tamed}
Let $M$ be a complete $G$-manifold. A {\em taming map} is a $G$-equivariant map $\v:M\to\grg$, such that the vector field $v$ on $M$, defined by \refe{v}, does not vanish anywhere outside of a cocompact subset of $M$. If $\v$ is a taming map, we refer to the pair $(M,\v)$ as a {\em tamed $G$-manifold}.

If, in addition, $\E$ is a $G$-equivariant $\ZZ_2$-graded self-adjoint Clifford module over $M$, we refer to the pair $(\E,\v)$ as a {\em tamed Clifford module} over $M$.
\end{Def}
The index we are going to define depends on the (equivalence class) of $\v$.

\begin{Rem}\label{R:family of metrics}
Suppose $M$ is a symplectic manifold and that the action of $G$ on $M$ is Hamiltonian with moment map $\mu:M\to \grg^*$. Following Mathai and Hochs \cite{MathaiHochs} we introduce a family of scalar products $\<\cdot,\cdot\>_x$ ($x\in M$) on $\grg^*$ which is $G$-invariant in the sense that
\[
	\<\Ad_g(\xi),\Ad_g(\eta)\>_{g\cdot x} \ = \ \<\xi,\eta\>_x, \qquad \text{for all}\quad x\in M, \ \xi,\eta\in \grg^*.
\]
For $x\in M$ set $\calH(x):=\<\mu(x),\mu(x)\>_x$ and let $\v(x)\in \grg$ denote the dual of $\mu(x)$ with respect to the scalar product $\<\cdot,\cdot\>_x$. Then the vector field \eqref{E:v} is equal to one half of the vector field $X_1^\calH$ introduced in \cite[\S2]{MathaiHochs}.  Thus the construction of the index below generalizes the construction of Mathai and Hochs to manifolds which are not symplectic. 
\end{Rem}

%---------------------------------------------
\subsection{A rescaling of $v$} \label{SS:rescaling}
Our definition of the index uses certain rescaling of the vector field  $v$. By this we mean the product $f(x)v(x)$, where $f:M\to[0,\infty)$ is a smooth positive function. Roughly speaking, we demand that $f(x)v(x)$ tends to infinity ``fast enough" when $q(x)\subset M/G$ tends to infinity. The precise conditions we impose on $f$ are quite technical, cf. \refd{admissible}. Luckily, our index turns out to be independent of the concrete choice of $f$. It is important, however, to know that at least one admissible function exists. This is guaranteed by \refl{rescaling} below.

We need to introduce some additional notations.

For a vector $\u\in\grg$, we denote by $\L^\E_\u$ the
infinitesimal action of $\u$ on $\Gam(M,\E)$ induced by the action
on $G$ on $\E$. Let $\n_{u}^\E:\Gam(M,\E)\to\Gam(M,\E)$ denote the covariant derivative along the vector field $u$ induced by $\u$. The difference between those two operators is a bundle map, which we denote by
\begin{equation}\label{E:mu}
    \mu^\E(\u) \ := \ \n^\E_{u}-\L^\E_\u \ \in \ \End \E.
\end{equation}

We will use the same notation $|\cdot|$ for  the norms on the
bundles $TM, T^*M, \E$.  Let $\End(TM)$ and $\End(\E)$ denote the
bundles  of endomorphisms of $TM$ and $\E$, respectively. We will
denote by $\|\cdot\|$ the norms on these bundles  induced by
$|\cdot|$. Let $\chi$ be a cutoff function as in Section~\ref{SS:cutoff}. To simplify the notation, set
\begin{equation}\label{E:nu}
    \nu=|\v|+\|\nLC v\|+\|\mu^\E(\v)\|+|v|+1.
\end{equation}
\begin{Def}\label{D:toinfty}
We say that a function $h:M\to [0,\infty)$ tends to infinity as $q(x)\to\infty$  and write $\lim_{q(x)\to\infty}h(x)= \infty$ if for arbitrary large number $R>0$ there exists a compact set $K\in M/G$ such that for all $x\not\in q^{-1}(K)$ we have $h(x)>R$.
\end{Def}

\begin{Def}\label{D:admissible}
We say that a smooth $G$-invariant function $f:M\to[0,\infty)$ on
a tamed $G$-manifold $(M,\v)$  is {\em admissible} for  $(\E,\v,\n^\E)$ if
\begin{equation}\label{E:limv}
        \lim_{q(x)\to\infty}\, \frac{f^2|v|^2}{
          |df||v|+f\nu+1 }  \ = \ \infty.
\end{equation}
\end{Def}
\lem{rescaling}
Let $(\E,\v)$ be a tamed Clifford module and let $\n^\E$ be a
$G$-invariant Clifford connection on $\E$. Let $\chi:M\to [0,\infty)$ be a cutoff function. Then  there exists an
admissible function $f$ for  $(\E,\v,\n^\E,\chi)$.
\elem
\begin{proof}
The proof for the case when the group $G$ is compact given in \cite[\S8]{Br-index} does not use the compactness of $G$ and works in our new situation.
\end{proof}
%----------------------------------------------
\subsection{Index on non-compact manifolds}\Label{SS:noncomind}
We use the Riemannian metric on $M$, to identify the tangent and the cotangent bundles of $M$. In particular, we consider $v$ as a
section of \/ $T^*M$.

Let $f$ be an admissible function and let $D_\chi$ be the operator defined in \eqref{E:Dphi}. Consider the {\em deformed
Dirac operator}
\begin{equation}\label{E:Dv}
        D_{\chi,fv} \ = \ D_\chi \ + \ {\i}c(fv):\, \chi\,\Gtc^G\ \to \ \chi\, \Gtc^G.
\end{equation}
This operator is essentially self-adjoint, cf. the remark on page 411 of \cite{Chernoff73}.

Our first result is the following analogue of Theorem~2.9 from \cite{Br-index}
\th{finite}
Suppose $f$ is an admissible function.  Then

1. \ The kernel of the deformed Dirac operator $D_{\chi,fv}$ has finite dimension. 

2. \ The index
\begin{equation}\label{E:indexDphifv}
	\ind_G D_{\chi,fv} \ := \ \dim\Ker D_{\chi,fv}^+ \ - \ \dim \Ker D_{\chi,fv}^-
\end{equation}	
 is independent of the choices of the cutoff function $\chi$, the  admissible function $f$, and the $G$-invariant Clifford connection $\n^\E$ on $\E$.
\eth
The proof of the first part of the theorem is given in \refs{prfinite}. The second part of the theorem is proven in  \refss{grch} as an immediate consequence of \reft{cobordinv} about cobordism invariance of the index.

We refer to the pair $(D,\v)$ as a {\em tamed Dirac operator}. The above theorem allows us to defined the index of a tamed Dirac operator $\ind_G(D,\v) := \ind_G(D_{\chi,fv})$.
Since $\ind_G(D,\v)$ is independent of the choice of the connection on $\E$, it is an invariant of the tamed Clifford module $(\E,\v)$. We set $\ind_G(\E,\v):=\ind_G(D,\v)$ and refer to it as   the {\em (analytic) index of a tamed Clifford module} $(\E,\v)$.

%-----------------------------------------
%----------------------------------------
\section{Cobordism invariance of the index}\Label{S:cobord}

In this section we adopt  the notion of cobordism between tamed Clifford modules and tamed Dirac operators introduced in  \cite{Br-index} to the case of a non-compact group $G$. We show that the index introduced in \refss{noncomind} is invariant under a cobordism.  We use this result to to prove Theorem~\ref{T:finite}.1.

%-------------------------
\subsection{Cobordism between tamed $G$-manifolds}\Label{SS:cobordM}
Note, first, that for cobordism to be meaningful  one must
make some compactness assumption. Otherwise, every manifold
is cobordant to the empty set via the noncompact cobordism
$M\times[0,1)$. Since our manifolds are non-compact
themselves, we can not demand cobordism to be compact.
Instead, we demand the cobordism to carry a taming map to
$\grg$.
\begin{Def}\label{D:cobordM}
A {\em cobordism \/} between tamed $G$-manifolds $(M_1,\v_1)$ and
$(M_2,\v_2)$ is a triple $(W,\v, \phi)$, where
\begin{enumerate}
  \item  $W$ is a complete Riemannian $G$-manifold with boundary;
  \item  $\v:W\to\grg$ is a smooth $G$-invariant map, such that
  the corresponding vector field $v$  does
  not vanish anywhere outside of a cocompact subset of $W$;
  \item  $\phi$ is a $G$-equivariant, metric preserving
  diffeomorphism between a neighborhood $U$ of the boundary $\d W$
  of $W$ and the disjoint union $\big(M_1\times[0,\eps)\big)\,
  \bigsqcup\, \big(M_2\times(-\eps,0]\big)$. We will refer to $U$
  as the {\em neck} and we will identify it with $\big(M_1\times[0,\eps)\big)\,
  \bigsqcup\, \big(M_2\times(-\eps,0]\big)$.
  \item the restriction of $\v\big(\phi^{-1}(x,t))$ to
  $M_1\times[0,\eps)$ (resp. to $M_2\times(-\eps,0]$)
  is equal to $\v_1(x)$ (resp. to $\v_2(x)$).
\end{enumerate}
\end{Def}

%------------------------------------------
\subsection{Cobordism between tamed Clifford modules}\Label{SS:cobordE}
If $M$ is a Riemannian $G$-manifold, then, for any interval $I\subset\RR$, the product $M\times I$ carries natural Riemannian
metric and $G$-action. Let $\pi:M\times I\to M, \ t:M\times{I}\to
I$ denote the natural projections. We refer to the pull-back
$\pi^*\E$ as a vector bundle {\em induced} by $\E$. We view $t$ as
a real valued function on $M$, and we denote by $dt$ its
differential.

%-------------
\begin{Def}\label{D:cobordE}
Let  $(M_1,\v_1)$ and $(M_2,\v_2)$ be tamed $G$-manifolds. Suppose
that each $M_i$, $i=1,2$, is endowed with a $G$-equivariant
self-adjoint Clifford module $\E_i=\E^+_i\oplus\E^-_i$. A {\em
cobordism \/} between the tamed Clifford modules  $(\E_i,\v_i)$,
$i=1,2$, is a cobordism $(W,\v,\phi)$ between $(M_i,\v_i)$
together with a pair $(\E_W,\psi)$, where
\begin{enumerate}
  \item $\E_W$ is a $G$-equivariant (non-graded) self-adjoint
  Clifford module over $W$;
  \item $\psi$ is a $G$-equivariant isometric isomorphism between
  the restriction of $\E_W$ to $U$ and the Clifford module induced
  on the neck $\big(M_1\times[0,\eps)\big)\bigsqcup
  \big(M_2\times(-\eps,0]\big)$ by $\E_i$.
  \item On the neck $U$ we have
    \(\displaystyle
        c(dt)|_{\psi^{-1}\E_i^\pm} \ = \ \pm \i.
    \)
\end{enumerate}
\end{Def}

\rem{M12-M}
Let $\E^{\text{op}}_1$ denote the Clifford module $\E_1$  with the
opposite grading, i.e., $\E^{\text{op}\pm}_1=\E_1^\mp$. Then,
$\ind_G(\E_1,\v_1)=  -\ind_G(\E^{\text{op}}_1,\v_1)$.

Consider the Clifford module $\E$ over the disjoint union
$M=M_1\sqcup{M_2}$ induced by the Clifford modules
$\E_1^{\text{op}}$ and $\E_2$. Let $\v:M\to\grg$ be the map such
that $\v|_{M_i}=\v_i$. A cobordism between $(\E_1,\v_1)$ and
$(\E_2,\v_2)$ may be viewed as  a cobordism between $(\E,\v)$ and
(the Clifford module over) the empty set.
\erem
\vskip0.2cm

One of the main results of this paper is the following theorem,
which asserts that the index is preserved by a cobordism.

\th{cobordinv}
   Suppose $(\E_1,\v_1)$ and $(\E_2,\v_2)$  are cobordant tamed Clifford
   modules. Let $D_1, D_2$ be Dirac operators associated to
   $G$-invariant Clifford connections on $\E_1$ and
   $\E_2$ and let $\chi_1,\chi_2$ be cutoff functions on $M_1$ and $M_2$. Then, for any admissible functions
   $f_1, f_2$ 
   \[
        \ind_G\big(\, D_{1,\chi_1}+{\i}c(f_1v_1)\, \big)
         \ = \ \ind_G\big(\, D_{2,\chi_2}+{\i}c(f_2v_2)\, \big).
   \]
\eth
The proof of the theorem is given in \refs{prcobordinv}.

%------------------------
\subsection{The definition of the analytic index of a tamed
Clifford module}\Label{SS:grch} \reft{cobordinv} implies, in
particular, that, if $(\E,\v)$ is a tamed Clifford module, then
the index $\ind_G(D_{\chi,fv})$ is independent of the choice of the
admissible function $f$, the cutoff function $\chi$,  and the Clifford connection on $\E$. This
proves part 2 of \reft{finite} and (cf. \refss{noncomind}) allows
us to define the {\em (analytic) index} of the tamed Clifford
module $(\E,\v)$
\[
    \ind_G(\E,\v) \ := \ \ind_G(D_{\chi,fv}), \qquad \text{$f$ is
    an admissible function}.
\]
\reft{cobordinv} can be reformulated now as
\th{cobordinva}
The indexes of cobordant tamed Clifford modules coincide.
\eth

%---------------------------------------------
\subsection{Index and zeros of $v$}\Label{SS:zerosv}
As a simple corollary of \reft{cobordinv}, we obtain the following 
\lem{zerosv} If the vector field $v(x)\not=0$ for all $x\in M$,
then $\ind_G(\E,\v)=0$.
\elem
\prf
Consider the product $W=M\times[0,\infty)$ and define the map
$\tv:W\to\grg$ by the formula: $\tv(x,t)=\v(x)$. Clearly,
$(W,\tv)$ is a cobordism between the tamed $G$-manifold $M$ and
the empty set. Let $\E_W$ be the lift of $\E$ to $W$. Define the
Clifford module structure $c:T^*W\to \End{\E_W}$ by the formula
\[
    c(x,a) e \ = \ c(x)e \ \pm \i a e, \qquad (x,a)\in T^*W\simeq
    T^*M\oplus\RR, \ e\in{\E_W^\pm}.
\]
Then $(\E_W,\tv)$ is a cobordism between $(\E,\v)$ and the
Clifford module over the empty set.
\eprf

%--------------------------------
\subsection{The excision property}\Label{SS:stDfv}
We will now amplify the above lemma and show that the index is
independent of the restriction of $(\E,\v)$ to a subset, where
$v\not=0$.

Let $(M_i,\v_i) \, i=1,2,$ be tamed $n$-dimensional $G$-manifolds.
Let $U$ be an open $n$-dimensional $G$-manifold. For each $i=1,2$,
let $\phi_i:U\to M_i$ be a smooth $G$-equivariant embedding. Set
$U_i=\phi_i(U)\subset M_i$. Assume that the boundary $\Sig_i=\d
U_i$ of $U_i$ is a smooth hypersurfaces in $M_i$. Assume also that
the vector field $v_i$ induced by $\v_i$ on $M_i$ does not vanish
anywhere on $M_i\backslash{U_i}$.

\lem{stabDv}
Let $(\E_1,\v_1)$,  $(\E_2,\v_2)$ be tamed Clifford modules over
$M_1$ and $M_2$, respectively. Suppose that the pull-backs
$\phi_i^*\E_i, \ i=1,2$ \/ are $G$-equivariantly isomorphic as
$\ZZ_2$-graded self-adjoint Clifford modules over $U$. Assume also
that $\v_1\circ\phi_1\equiv \v_2\circ\phi_2$. Then $(\E_1,\v_1)$
and $(\E_2,\v_2)$ are cobordant. In particular,
$\ind_G(\E_1,\v_1)= \ind_G(\E_2,\v_2)$.
\elem
\begin{proof}
An explicit cobordism between $(\E_1,\v_1)$ and $(\E_2,\v_2)$ is constructed in Section~12.2 of \cite{Br-index}.
\end{proof}

The following lemma is, in a sense, opposite to \refl{stabDv}. 

\lem{vcomp}
Let $\v_1,\v_2:M\to \grg$ be taming maps, which coincide out of a
compact subset of $M$. Then the tamed Clifford modules $(\E,\v_1)$
and $(\E,\v_2)$ are cobordant. In particular,
$\ind_G(\E,\v_1)= \ind_G(\E,\v_2)$.
\elem
\prf
The proof is a verbatim repetition of the proof of  Lemma~3.16 in \cite{Br-index}.
\eprf

%-------------------------------------------------------
%-------------------------------------------------------
\section{The gluing formula}\Label{S:gluing}

If we cut a tamed $G$-manifold along a $G$-invariant hypersurface
$\Sig$, we obtain a manifold with boundary. By rescaling the
metric near the boundary we may convert it to a complete manifold
without boundary, in fact, to a tamed $G$-manifold. In this
section, we show that the index is invariant under this type of
surgery. In particular, if $\Sig$ divides $M$ into two pieces
$M_1$ and $M_2$, we see that the index on $M$ is equal to the sum
of the indexes on $M_1$ and $M_2$. In other words, the index is
{\em additive}. 

%------------------------------
\subsection{The surgery}\Label{SS:surgery}
Let $(M,\v)$ be a tamed $G$-manifold. Suppose $\Sig\subset M$ is a
smooth $G$-invariant hypersurface in $M$. For simplicity, we
assume that $\Sig$ is cocompact. Assume also that the vector field
$v$ induced by $\v$ does not vanish anywhere on $\Sig$. Suppose
that $\E=\E^+\oplus\E^-$ is a $G$-equivariant $\ZZ_2$-graded
self-adjoint Clifford module over $M$. Denote by $\E_\Sig$ the
restriction of the $\ZZ_2$-graded Hermitian vector bundle $\E$ to
$M_\Sig:= M\backslash\Sig$.

Let $g^M$ denote the Riemannian metric on $M$. In Section~4.2 of \cite{Br-index} by a rescaling of $g^M$ we constructed a complete Riemannian metric $g^{M_\Sig}$ on $M_\Sig$ and a Clifford action $c_\Sig:T^*M_\Sig\to \End(\E_\Sig)$ compatible with this metric. It follows, from the cobordism invariance of the index (\refl{stabDv}), that  our index theory is essentially independent of the concrete choice of $g^{M_\Sig}$.

%%%%%%%
\th{gluing}
The tamed Clifford modules $(\E,\v)$ and $(\E_\Sig,\v_\Sig)$ are
cobordant. In particular,
\[
    \ind_G(\E,\v) \ = \ \ind_G(\E_\Sig,\v_\Sig).
\]
\eth
We refer to \reft{gluing} as a {\em gluing formula}, meaning that
$M$ is obtained from $M_\Sig$ by gluing along $\Sig$.

\begin{proof}
To prove the theorem it is enough to construct a cobordism  $W$
between $M$ and $M_\Sig$.

Consider the product $M\times[0,1]$, and the set
\[
    Z \ := \ \big\{\, (x,t)\in M\times[0,1]: \, t\le 1/3, x\in \Sig\,
    \big\}.
\]
Set $W:= (M\times[0,1])\backslash{Z}$. Then $W$ is a $G$-manifold,
whose boundary is diffeomorphic to the disjoint union of
$M\backslash\Sig\simeq (M\backslash\Sig)\times\{0\}$ and $M\simeq
M\times\{1\}$. To finish the proof we need to construct a complete Riemannian
metric $g^W$ on $W$, so that the condition (iii) of \refd{cobordM}
is satisfied.  For the case when the group $G$ is compact it is done in Section~13 of \cite{Br-index}. The same construction works without any changes for a non-compact $G$.
\end{proof}

%----------------------------------
\subsection{The additivity of the index}\Label{SS:add}
Suppose that $\Sig$ divides  $M$ into two open submanifolds $M_1$
and $M_2$, so that $M_\Sig=M_1\sqcup M_2$. The metric $g^{M_\Sig}$
induces complete $G$-invariant Riemannian metrics $g^{M_1},
g^{M_2}$ on $M_1$ and $M_2$, respectively. Let $\E_i, \v_i \
(i=1,2)$ denote the restrictions of the Clifford module $\E_\Sig$
and the taming map $\v_\Sig$ to $M_i$.  Then \reft{gluing} implies
the following
\cor{gluing}
$\displaystyle \ind_G(\E,\v) \ = \
        \ind_G(\E_1,\v_1) \ + \ \ind_G(\E_2,\v_2)$.
\ecor
Thus, we see that the index of non-compact manifolds is {\em
``additive"}.

%----------------------------------------------------------
%----------------------------------------------------------

\section{Functions on a  cobordism}\Label{S:prrescaling}

In this section we define the notions of  an  admissible and cutoff functions on a cobordism and prove
the existence of such functions. 

\subsection{A cobordism}\Label{SS:admiscob}
Let $(\E,\v)$ be a tamed Clifford module over a complete
$G$-manifold $M$. Let $(W,\v_W,\phi)$ be a cobordism between
$(M,\v)$ and the empty set, cf. \refd{cobordM}. In particular, $W$
is a complete $G$-manifold with boundary and $\phi$ is a
$G$-equivariant metric preserving diffeomorphism between a
neighborhood $U$ of $\d{W}\simeq M$ and the product
$M\times[0,\eps)$.

Let $\pi:M\times[0,\eps)\to M$ be the projection. A $G$-invariant
Clifford connection $\n^{\E}$ on $\E$ induces a connection
$\n^{\pi^*\E}$ on the pull-back $\pi^*\E$, such that
\begin{equation}\label{E:tnE}
    \n^{\pi^*\E}_{(u,a)} \ := \ \pi^*\n^{\E}_u \ + \
       a\frac{\d}{\d t}, \qquad
            (u,a)\in TM\times \RR\simeq T(M\times[0,\eps)).
\end{equation}
Let $(\E_W,\v_W,\psi)$ be a cobordism between $(\E,\v)$ and the
unique Clifford module over the empty set, cf. \refd{cobordE}. In
particular, $\psi:\E_W|_U\to \pi^*\E$  is a $G$-equivariant
isometry. Let $\n^{\E_W}$ be a $G$-invariant connection on $\E_W$,
such that
$\n^{\E_W}|_{\phi^{-1}(M\times[0,\eps/2))}=\psi^{-1}\circ\n^{\pi^*\E}\circ\psi$.

%----------
\subsection{An admissible function on a cobordism}\label{SS:admiscobord}

%--------
\begin{Def}\label{D:admisb}
A smooth $G$-invariant function $f:W\to[0,\infty)$ is an admissible function for $(\E_W,\v_W,\n^{\E_W},\chi)$, if
it satisfies \refe{limv} and there exists a function $h:M\to[0,\infty)$ such that $f\big(\phi^{-1}(y,t)\big)= h(y)$ for all $y\in M, t\in[0,\eps/2)$.
\end{Def}
%-------
\lem{admisb}
Suppose $h$ is an admissible function for $(\E_{M},\v,\n^{\E})$.
Then there exists an admissible function $f$ on
$(\E_W,\v_W,\n^{\E_W})$ such that the restriction $f|_{M}=h$.
\elem
\prf
For the case when $G$ is compact the proof is given in Section~8 of \cite{Br-index}. The proof does not use the compactness of $G$ and extends to non-compact case without any changes. 
\eprf
%----------
\subsection{A cutoff function on a cobordism}\label{SS:cutoffcobord}
%----------------------
\begin{Def}\label{D:cutoffW}
A smooth function $\chi:W\to [0,\infty)$ is called a {\em cutoff function on cobordism} if its support intersects all $G$-orbits in compact sets, it satisfies  \eqref{E:cutoff}, and there exists a cutoff function $\eta:M\to [0,\infty)$ such that  $\chi\big(\phi^{-1}(y,t)\big)= \eta(y)$ for all $y\in M, t\in[0,\eps/2)$.
\end{Def}

%-----
\lem{cutoffb}
Suppose $\eta$ is a cutoff function on $M$. Then there exists a cutoff function $\chi$ on $W$, such that $\chi|_M=\eta$.
\elem
\begin{proof}

Recall that we are given  a neighborhood $U$ of $\partial{W}$ and a diffeomorphism $\phi: U\to M\times{}I$. We write $\phi(x)=(y,t)$ where $y\in M$ and $t\in [0,\eps)$.  

Let $\alp$ and $\bet$ be smooth functions $[0,\infty)\to [0,1]$ such that 
\[
	\alp(t)\ = \ \begin{cases} 0, \quad\text{for}\ \ t&<\frac{\eps}3, \\ 1, \quad\text{for} \ \  t&>\frac{2\eps}3. \end{cases}
\]
and 
\begin{equation}\label{E:alp2+bet2}
	\alp^2+\bet^2 \ \equiv\ 1.
\end{equation} 
Let  $\chi_1$ be any cutoff function on $W$ and define a new function $\chi_2$ on $W$ by
\[
	\chi_2(x) \ = \ \begin{cases}
	   \chi_1(x), \quad&\text{for}\ \ x\not\in U,\\
	   \alp(t)\chi_1(x) + \bet(t)\eta(y), \quad&\text{for}\ \ x=(y,t)\in U\simeq M\times I.
	  \end{cases}
\]
Then the support of $\chi_2$ intersects the  orbits of $G$ in compact sets and 
\[
	\chi_2^2 \ = \ \alp^2\chi_1^2+\bet^2\eta^2+2\alp\bet\chi_1\eta \ = \ \alp^2\chi_1^2+\bet^2\eta^2 + \psi,
\]
where we set $\psi:= 2\alp\bet\chi_1\eta$. Hence, form \eqref{E:alp2+bet2} and the definition of the cutoff function \eqref{E:cutoff}  we obtain 
\[
	\int_G\, \chi_2^2(g\cdot x) \,dg \ = \ 1\ +\ \int_G\,\psi(g\cdot x)\,dg.
\]
It follows that the function 
\[
	\chi(x)\ := \ \frac{\chi_2(x)}{ \big(\,1\ +\ \int_G\psi(g\cdot x)dg\,\big)^{1/2}}
\]
is a cutoff function whose restriction to $\partial{W}$ is equal to $\eta$.
\end{proof}

%------------------------------------------------------
%------------------------------------------------------

\section{Proof of \reft{finite}.1}\Label{S:prfinite}

%-----------------------------------------------
\subsection{The operator $D_{fv}^2$}\label{SS:Du2til}
Consider the operator 
\[
	D_{fv}\ : = \ D+\sqrt{-1}c(fv).
\]
It follows from \eqref{E:Dphi} that 
\begin{equation}\label{E:Dchiphi}\notag
	D_{\chi,fv}(\chi s)\ := \ \chi\,D_{fv}s, \qquad s\in \Gtc,
\end{equation}
and, hence, 
\begin{equation}\label{E:D2chiphi}
	D_{\chi,fv}^2(\chi s)\ := \ \chi\,D_{fv}^2s, \qquad s\in \Gtc,
\end{equation}

%-----------------------------------------------
\subsection{Calculation of $D_{\chi,fv}^2$}\Label{SS:Du2}
Let $f$ be an admissible function and set $u=fv$. Consider the
operator
\begin{equation}\label{E:Av}
        A_u \ = \ \sum c(e_i)\, c(\nLC_{e_i}u):\, \E \ \to \ \E,
\end{equation}
where $\oe=\{e_1\ldots e_n\}$ is an orthonormal frame of
$TM\simeq T^*M$ and $\nLC$ is the Levi-Civita connection on
$TM$. One easily checks that $A_u$ is independent of the
choice of $\oe$ (it follows, also, from equation \eqref{E:D2} below).

Recall that we denote by $\calL_\bfu$ the infinitesimal action of $\bfu$ on $\Gam(M,\E)$ induced by the action of $G$ on $\calE$. Then the restriction of $\calL_\bfu$ to $\Gtc^G$ is equal to 0. Combining Lemma~9.2 of \cite{Br-index} with \eqref{E:mu} we obtain
\begin{equation}\label{E:D2}
    D_{u}^2 \ = \ D^2 \ + \ |u|^2 \ + \ {\i}A_u \ - \ 2{\i}\mu^\E(\bfu):\, \Gtc^G\to \Gtc^G.
\end{equation}

Since both $A_u$ and $\mu^\E(\bfu)$ commute with multiplication by $\chi$, we conclude from \eqref{E:Dphi}, \eqref{E:D2chiphi}, and \eqref{E:D2} that 
\begin{equation}\label{E:Dchi2}
    D_{\chi,u}^2 \ = \ D_\chi^2 \ + \ |u|^2 \ + \ {\i}A_u \ - \ 2{\i}\mu^\E(\bfu):\, 
    \chi\Gtc^G\to \chi\Gtc^G.
\end{equation}

%---------------------------
\subsection{Proof of \reft{finite}}\Label{SS:prfinite}
It is shown in  Section~9.3 of \cite{Br-index}  that  there exists a real valued function $r(x)$ on $M$ such that
$\lim_{q(x)\to\infty}\, r(x) \ = \ +\infty$  and
\[
	|u|^2 \ + \ {\i}A_u \ - \ 2{\i}\mu^\E(\bfu) \ \ge \ r(x).
\]
It follows now from \eqref{E:Dchi2} that
\begin{equation}\label{E:Du2a}
     D_{\chi,u}^2\ \ge \ D_\chi^2 \ + r(x).
\end{equation}
By Proposition~3.9 of \cite{MathaiHochs} the operator $D_{\chi,u}$ is Fredholm. 
\hfill$\square$

%-------------------------------------------------------
%--------------------------------------------------------
\section{Proof of \reft{cobordinv}}\Label{S:prcobordinv}

The proof is a minor modification of  the proof of Theorem~3.7 in \cite{Br-index}. We only sketch it here, stressing the couple of places where it defers from \cite{Br-index}. 

By \refr{M12-M}, it is enough to show that, if $(\E,\v)$ is
cobordant to (the Clifford module over) the empty set, then
$\ind_G(D_{\chi,fv})=0$ for any admissible function $f$.

Let $(W,\E_W,\v)$ be a cobordism between the empty set and
$(\E,\v)$ (slightly abusing the notation, we denote by the same
letter $\v$ the taming maps on $W$ and $M$).

In Section~\ref{S:prrescaling} we showed that there are exist a cutoff function and an admissible function on $W$ whose restriction to $M$ are equal to $\chi$ and $f$ respectively. By a slight abuse of notation, we denote these functions by the same letters $\chi$ and  $f$.

Let $\tilW$ be the manifold obtained from $W$ by attaching a
cylinder to the boundary, i.e.,
\[
    \tilW \ = \ W \, \sqcup \,
                \big(\, M\times(0,\infty)\, \big).
\]
The action of $G$, the Riemannian metric,  the map $\v$, the
functions $\chi, \ f$ and the Clifford bundle $\E_W$ extend naturally from
$W$ to $\tilW$.

Let us consider two anti-commuting actions (left and right action) of  $\RR$ on the exterior algebra
$\Lam^\b\CC=\Lam^0\CC\oplus\Lam^1\CC$, given by the formulas
\eq{ClonR}
    c_L(t)\, \ome \ = \ t\wedge\ome \ - \ \iot_t\ome; \qquad
    c_R(t)\, \ome \ = \ t\wedge\ome \ + \ \iot_t\ome.
\end{equation}
Note, that $c_L(t)^2=-t^2$, while $c_R(t)^2=t^2$. 

Set $\wE= \E\otimes\Lam^\b\CC$ and define the grading and the Clifford action $\tilc:T^*\tilW\to\End\wE$ on $\wE$   by the formulas
\[
   \wE^+ \ := \ \E_W\otimes\Lam^0\CC; \quad  \wE^- \ := \
   \E_W\otimes\Lam^1\CC; \qquad
   \tilc(v) \ := \ \i c(v)\otimes c_L(1) \quad (v\in T^*\tilW).
\]

Let $\tilD$ be a Dirac operator on $\wE$ and consider the operator
 \[
 	\tilD_{\chi,fv}\ :=\ \tilD_\chi+\i{}\tilc(fv).
\]
Note that $v$ might vanish somewhere near
infinity on the cylindrical end of $\tilW$. In particular, the index of $\tilD_{\chi,fv}$ is
not defined in general.

Let $p:\tilW\to\RR$ be a map, whose restriction to
$M\times(1,\infty)$ is the projection on the second factor,
and such that $p(W)=0$. For every $a\in\RR$, consider the
operator 
\[
	\bfD_a\ :=\ \tilD_{\chi,fv}-1\otimes c_R((p(t)-a)).
\] 
By Lemma~10.4 of \cite{Br-index} we have 
\begin{equation}\label{E:bfD2}
    \bfD_a^2 \ = \ \tilD_{\chi,fv}^2 - B + |p(x)-a|^2,
\end{equation}
where $B:\wE\to\wE$ is a uniformly bounded bundle map.

It follows now  from Proposition~3.9 of \cite{MathaiHochs}  that the operator $\bfD_a$ is Fredholm. 

%-----
\prop{depofa}
\(\displaystyle\ind_G(\bfD_a)=0 \) for all $a\in\RR$.
\eprop
\prf
Since
\[
    \bfD_a-\bfD_b\ = \ 
    c_R(b-a):\, \chi\Gam_{tc}^\infty(\tilW,\tilE)\ \to \ \chi\Gam_{tc}^\infty(\tilW,\tilE)
\]
is bounded operator depending continuously on $a,b\in\RR$, the index $\ind_G(\bfD_a)$ is independent of $a$. Therefore, it is enough to prove the proposition for one particular value of $a$. 

Let $\|B(x)\|, \, x\in\tilW$, denote the norm of the bundle map $B_x:\wE_x\to \wE_x$ and let $\|B\|_\infty= \sup_{x\in\tilW}\|B(x\|$. Choose $a\ll0$ such that $a^2>\|B\|_\infty$. Since $p(x)\ge0$ we  have $|p(x)-a|^2\ge a^2> \|B\|_\infty$. Thus it follows from \eqref{E:bfD2} that  $\bfD_a^2>0$, so that $\Ker\bfD^2_a=0$. Hence,
$\ind_G(\bfD^2_a)=0$.
\eprf

\reft{cobordinv} follows now from \refp{depofa} and the following
\th{a-a}
\(\displaystyle \ind_G(\bfD_a) \ = \
         \ind_G(D_{\chi,fv}) \).
\eth
\begin{proof}
For compact group $G$ the proof is given in Section~11 of \cite{Br-index}. This proof works without any changes in our current situation. 
\end{proof}

%------------------------------------------------
%-----------------------------------------------
\def\cprime{$'$} \def\cprime{$'$} \newcommand{\noop}[1]{} \def\cprime{$'$}
\providecommand{\bysame}{\leavevmode\hbox to3em{\hrulefill}\thinspace}
\providecommand{\MR}{\relax\ifhmode\unskip\space\fi MR }
\providecommand{\MRhref}[2]{%
  \href{http://www.ams.org/mathscinet-getitem?mr=#1}{#2}
}
\providecommand{\href}[2]{#2}


\begin{thebibliography}{10}

\bibitem{BeGeVe}
N.~Berline, E.~Getzler, and M.~Vergne, \emph{Heat kernels and {Dirac}
  operators}, Springer-Verlag, 1992.

\bibitem{Bourbaki-integration7-8}
N.~Bourbaki, \emph{\'{E}l\'ements de math\'ematique. {F}ascicule {XXIX}.
  {L}ivre {VI}: {I}nt\'egration. {C}hapitre 7: {M}esure de {H}aar. {C}hapitre
  8: {C}onvolution et repr\'esentations}, Actualit\'es Scientifiques et
  Industrielles, No. 1306, Hermann, Paris, 1963.

\bibitem{Br-index}
M.~Braverman, \emph{Index theorem for equivariant {D}irac operators on
  noncompact manifolds}, $K$-Theory \textbf{27} (2002), no.~1, 61--101.

\bibitem{BrBackground}
\bysame, \emph{Background cohomology of a non-compact {K\"ahler}
  ${G}$-manifold}, Trans. Amer. Math. Soc. \textbf{367} (2014).

\bibitem{BrCano}
M.~Braverman and L.~Cano, \emph{Index theory for non-compact ${G}$-manifolds},
  Geometric, Algebraic and Topological Methods for Quantum Field Theory, World
  Scientific Publishing Company, 2013, pp.~60--94.

\bibitem{Chernoff73}
P.~Chernoff, \emph{Essential self-adjointness of powers of generators of
  hyperbolic equations}, J. Functional Analysis \textbf{12} (1973), 401--414.

\bibitem{GGK96}
V.~L. Ginzburg, V.~Guillemin, and Y.~Karshon, \emph{Cobordism theory and
  localization formulas for {Hamiltonian} group actions}, Internat. Math. Res.
  Notices \textbf{5} (1996), 221--234.

\bibitem{GromLaw83}
M.~Gromov and B.~Lawson, \emph{Positive scalar curvature and the {D}irac
  operator on complete {R}iemannian manifolds}, Inst. Hautes \'Etudes Sci.
  Publ. Math. (1983), no.~58, 295--408.

\bibitem{GGK-book}
V.~Guillemin, V.~Ginzburg, and Y.~Karshon, \emph{Moment maps, cobordisms, and
  {H}amiltonian group actions}, Mathematical Surveys and Monographs, vol.~98,
  American Mathematical Society, Providence, RI, 2002, Appendix J by Maxim
  Braverman.

\bibitem{MathaiHochs}
P.~Hochs and V.~Mathai, \emph{Geometric quantization and families of inner products}, 50 pages, [arXiv:1309.6760].

\bibitem{Karshon98}
Y.~Karshon, \emph{Momentum maps and non-compact cobordisms}, J. Diff. Geom.
  \textbf{49} (1998), 183--201.

\bibitem{MaZhang-noncompact}
X.~Ma and W.~Zhang, \emph{Geometric quantization for proper moment maps}, 2008,
  arXiv:0812.3989.

\bibitem{MaZhang_TrIndex12}
\bysame, \emph{Transversal index and ${L}_2$-index for manifolds with
  boundary}, Progress in Mathematics, vol. 297, Birkh\"auser Boston Inc.,
  Boston, MA, 2012, pp.~299--316.

\bibitem{MathaiZhang10}
V.~Mathai and W.~Zhang, \emph{Geometric quantization for proper actions}, Adv.
  Math. \textbf{225} (2010), 1224--1247, [arXiv:0806.3138]

\bibitem{Paradan03}
P.-{\'E}. Paradan, \emph{{$\operatorname{Spin}\sp c$}-quantization and the
  {$K$}-multiplicities of the discrete series}, Ann. Sci. \'Ecole Norm. Sup.
  (4) \textbf{36} (2003), 805--845.

\bibitem{Vergne07}
M.~Vergne, \emph{Applications of equivariant cohomology}, International
  {C}ongress of {M}athematicians. {V}ol. {I}, Eur. Math. Soc., Z\"urich, 2007,
  pp.~635--664.

\end{thebibliography}
\end{document}